\documentclass[12pt]{amsart}
\usepackage{thmtools,thm-restate} 
\usepackage[margin=1.1in]{geometry}
\usepackage{url}
\usepackage{todonotes}

\usepackage{amssymb, amsmath, amscd, amsthm, enumitem, algpseudocode, algorithm, mathtools, comment}
\usepackage{amsfonts}
\usepackage[all]{xy}
\usepackage{graphicx}
\usepackage[normalem]{ulem}
\usepackage{hyperref}
\usepackage{faktor}

\newtheorem{thm}{Theorem}
\newtheorem{cor}{Corollary}[section]
\newtheorem{lem}{Lemma}[section]
\newtheorem{prop}{Proposition}[section]
\newtheorem*{Lemma 5.3}{Lemma 5.3}
\theoremstyle{definition}
\newtheorem{defn}{Definition}[section]
\theoremstyle{remark}
\newtheorem{rem}{Remark}[section]

\numberwithin{equation}{section}
\newtheorem{ex}{Example}[section]
\newtheorem{ass}{Assumption}[section]
\newtheorem*{thm*}{Lemma 5.8}

\newcommand{\R}{\mathbb{R}}
\newcommand{\N}{\mathbb{N}}
\newcommand{\Z}{\mathbb{Z}}
\newcommand{\K}{\mathbb{K}}
\newcommand{\p}{\varphi}

\newcommand{\s}{\psi}
\newcommand{\eps}{{\varepsilon}}

\newcommand{\hcap}{\; \hat{\cap} \;}

\newcommand{\Dgm}[1]{
    \textup{Dgm}\left(#1\right)
}
\newcommand{\dgm}[3]{
	\Dgm{#1^*_{(#2,#3)}}
}
\newcommand{\semip}{
	\bar{\Delta}^*
} 
\newcommand{\dmatch}[2]{
    d_{\textup{B}}\left(#1,#2\right)
}
\newcommand{\Dmatch}[2]{
    D_{\textup{match}}\left(#1,#2\right)
}
\newcommand{\norm}[2]{
    \|#1\|_{#2}
}
\newcommand{\spe}[2]{
    \textup{Sp}\left(#1,#2\right)
}

\newcommand{\Rinf}[0]{\mathbb{R}_{{-}\!\!{\bullet}}
}

\title[Matching distance via the extended Pareto grid]{Matching distance via the extended Pareto grid }

\author{Patrizio Frosini}
\address{Department of Mathematics, University of Bologna, Italy}
\email{patrizio.frosini@unibo.it}
\author{Eloy M\'osig Garc\'\i a}
\address{Department of Mathematics, University of Pisa, Italy}
\email{eloy.mosig@phd.unipi.it}
\author{Nicola Quercioli}
\address{DEI \& WiLab-National Laboratory for Wireless Communications, CNIT, University of Bologna, Italy}
\email{nicola.quercioli2@unibo.it}
\author{Francesca Tombari}\address{Max Planck Institute for Mathematics in the Sciences, Germany}
\email{francesca.tombari@mis.mpg.de}

\begin{document}
\maketitle

\begin{abstract}
One of the most animated themes of multidimensional persistence is the comparison of invariants.
The matching distance between persistent Betti numbers functions (or rank invariants) is among the most studied metrics in this context, particularly in 2-parameter persistence. 
The main reason for this interest is that, in the 2-parameter case, the foliation method allows for an effective computation of the matching distance, based on filtering the space along lines of positive slope.
Our work provides a qualitative analysis, based on a construction called extended Pareto grid, of the filtering lines that contribute to the computation of the matching distance. 
Under certain genericity assumptions, we show that these lines must either have slope 1 or belong to a collection of special lines associated with discontinuity phenomena.

\end{abstract}

\section{Introduction}
\label{sec:introduction}
One of the main tasks of Topological Data Analysis, in general, and persistent homology, in particular, is data comparison. 
It is a well-known observation that filtering a topological space by its sublevel sets gives finer homological information about the space than that contained in the total homology. 
In fact, taking the homology of a filtered space also captures some aspects of the geometry of the space, not only its global topology. 
Data comparison, in this sense, is done at the level of the topological and geometrical information inferred from the space and encoded in some compact and discriminative invariants. 

The matching distance, introduced in \cite{Biasotti2008} in homological degree 0 and in \cite{pers_Betti_stable} in general degree, is our candidate for such a task. 
This is a metric between persistent Betti numbers functions, or rank invariants, of continuous filtering functions. 
As such, it is defined in the context of $1$-parameter persistent homology, where it coincides with the bottleneck distance, but also in the multiparameter case. 
Of particular interest for us is its definition in $2$-parameter persistence.
There are two main reasons for this, and the first one is computational. 
In recent years,  strategies and algorithms to compute and approximate the matching distance efficiently have been object of research (see for example \cite{comp_disc,comp2,comp_asymp, comp_exact} and \cite{comp1, efficient-approx}). 
However, one of the most interesting aspects for this work is the availability of a tool in $2$-parameter persistence, called the \emph{extended Pareto grid}. 
The extended Pareto grid is a union of closed arcs and half-lines in the plane associated with a class of smooth $\R^2$-valued function on a Riemannian manifold. 
This concept was introduced in \cite{coher_match} as a way to study biparameter persistence, as it has been proven a powerful tool for tracking homological features (see also \cite{EPG_pers}).
Its origins, however, lie in Morse theory \cite{wan}.

The $2$-parameter matching distance between two functions is defined as a supremum over all filtering lines of positive slope on the plane along which the underlying space is filtered. 
Each one of these $1$-parameter filtrations gives a persistence diagram of the restriction of the functions we want to compare, between which it is possible to compute the bottleneck distance. 
In \cite{slope1}, the authors give a qualitative characterisation of the filtering lines where the supremum defining the bottleneck distance can be attained. 
Their results assert that such a supremum is attained at a vertical, horizontal or slope 1 line, or at a class of lines associated with \emph{special values} depending on the extended Pareto grids of the two functions. 
Our work builds on the shoulder of \cite{slope1} and states that, under certain assumptions on the filtering functions and their extended Pareto grids, the matching distance is attained either at a line of slope 1, or at a subset of lines associated with special values that is, in general, smaller than the one in \cite{slope1} (see Theorem \ref{main}). 
Furthermore, we conjecture that such subset of special values is generically finite.

This work is not the first attempt to study the $2$-parameter matching distance from a computational point of view. 
However, together with~\cite{coher_match, slope1}, it is the first one, to our knowledge,
that considers the matching distance between persistence modules that are not finitely generated.
A similar approach to ours is used in~\cite{assivbary, EPG_pers} to study the persistence modules associated with 2-parameter filtrations. 
To our knowledge, it is still unknown if such persistence modules, despite not being finitely generated, are generically finitely encoded in the sense of \cite{miller2020homological}.


\section{Mathematical setting}
\label{sec:math-setting}
In this section we describe 1-parameter and 2-parameter persistent homology of a filtering function, together with their invariants and extended metrics between such invariants.
We also introduce a key construction for this article, called the extended Pareto grid. 
Most of the content of this section is not novel and can be found in the literature (see, for example,~\cite{pers_Betti_stable}).

\subsection{$1$-parameter Persistent Homology}

Let $X$ be a compact topological space and $\p \colon X \to \R$ a continuous function.
The sublevel set of $\p$ at $u$, $\{(x\in X\mid \varphi(x)\le u)\}$ is denoted by $(\p\le u)$.
The function $\p$ induces a filtration on $X$ given by its sublevel sets, $(\p\le u)\subseteq(\p\le v)$, for $u\le v$.
For any $k \in \Z$ and a field $\K$, we can consider the Čech homology functor $H_k\colon \text{Top}\to \text{vect}_{\K}$, which associates the sublevel set $(\p\le v)$ to the $\K$-vector space $H_k(\p\le v)$, and the inclusion $\iota_{u,v} \colon (\p\le u)\subseteq (\p\le v)$ to the linear map $\iota^*_{u,v} \colon H_k(\p\le u) \to H_k(\p\le v)$. 
From now on, whenever a continuous function from a compact set is considered, we will require the following to hold:
\begin{ass}\label{jdshdsds}
$H_k(\p\le u)$ is a finitely generated $\K$-vector space for every degree $k \in \Z$ and for every $u \in \R$.
\end{ass}

The \textbf{$k$-th persistent homology group} of $\p$ at $(u,v)$ is defined as $H_k^\p(u,v)= \textup{Im }\iota^*_{u,v} \subseteq H_k(\p\le v)$. 
The function $\beta^\p_k(u,v) = \textup{dim } H_k^\p(u,v)$, defined for $u\le v$, is called \textbf{$k$-th persistent Betti numbers function}, or PBNF, of $\p$ at $(u,v)$.
Note that $\beta^\p_k(u,v)$ is non-decreasing in the first variable and non-increasing in the second (see~\cite{pers_Betti_stable}).

Let $\Delta^+$ denote the half-plane $\{(u,v) \in \R^2 \mid u < v \}$ and $\Delta^* $ the union $ \Delta^+\cup\{(u,\infty)\mid u\in\R\}$.
Moreover, set $\Delta$ to be the diagonal $\{(u,v) \in \R^2 \mid u = v \}$. 
The set $\bar\Delta^*$ is defined as the quotient $(\Delta^* \cup\Delta)/_\sim$, where $(u,v)\sim (u',v')$ if, and only if, $u=v$ and $u'=v'$, or $u=u'$ and $v=v'$.
We denote by $\Delta$ the equivalence class $[(u,u)]_\sim$, for any $u$. 
With a little abuse of notation, we also denote by $(u,v)$ the equivalence class with one element $(u,v)\in \Delta^\ast$.
We endow $\bar\Delta^*$ with the metric $d$ defined as
\begin{equation}\label{opcklfajkclòkweaa}
d(p,q)=\begin{cases}
             C(u,u',v,v') & \textup{if } p = (u,v),\; q = (u',v') \in \Delta^+, \\
            \vert u - u' \vert & \textup{if } p = (u,\infty),\; q = (u',\infty), \\
            \frac{v-u}{2} & \textup{if } p = (u,v) \in \Delta^+, \; q = \Delta, \\
            \frac{v'-u'}{2} & \textup{if } q = (u',v') \in \Delta^+, \; p = \Delta, \\
            0 & \textup{if } p = q = \Delta, \\
            \infty & \textup{otherwise},
            \end{cases}
\end{equation}
for every $p,q$ in $\bar\Delta^*$, and where $C(u,u',v,v')$ denotes the value $\min \{ \max \{ \vert u - u' \vert, \vert v - v' \vert\},\allowbreak \max \{ \frac{v-u}{2}, \frac{v'-u'}{2} \}\}$. 
Notice, in particular, that every point $(u, \infty)\in \Delta^\ast\setminus \Delta^+$ is at infinite distance from any point $(u',v')\in \Delta^+$ and from $\Delta$.

Let us consider $\eps>0$ and $u < v < \infty$. 
We define the number 
\[
\mu_{\eps, k}^\p (u,v) = \beta^\p_k(u + \eps,v-\eps) - \beta^\p_k(u - \eps,v-\eps) + \beta^\p_k(u-\eps,v+\eps) - \beta^\p_k(u + \eps,v+\eps)
\]
and call it the \textbf{total multiplicity} of the $\eps$-box centered at $(u,v)$.
Then the \textbf{k-th persistence diagram} of the continuous function $\p \colon X \to \R$ is the multiset
$\mathrm{Dgm}_k(\p)   = (\semip, m)$,
\[
m(u,v) = 
\begin{cases}
\mu(u,v) = \textup{lim}_{\eps \to 0^+} \; \mu_{\eps, k}^\p (u,v) &\textup{if } u < v < \infty\\
\nu(u,v) = \textup{lim}_{\eps \to 0^+} \; \beta_k^\p (u + \eps,v) - \beta_k^\p (u - \eps,v) &\textup{if } u < v = \infty\\
\infty &\textup{if } (u,v) = \Delta.
\end{cases}
\]
The number $m(u,v)$ is called the \textbf{multiplicity} of the point $(u,v)$.
We often just write $\mathrm{Dgm}(\p)$, if a statement holds for any degree of homology.
We also refer to the points of the persistence diagram with positive multiplicity as \textbf{cornerpoints}. 
Those cornerpoints $(u,v)$ in $\text{Dgm}(\p)$ with finite coordinates are called \textbf{proper cornerpoints}; while those of the form $(u, \infty)$ are called \textbf{improper} or \textbf{essential cornerpoints}. The point $\Delta$ is the trivial cornerpoint.
With a further abuse of notation, we also denote by $\text{Dgm}(\p)$ the collection of points of $\semip$ with positive multiplicity. 

The space of persistence diagrams can be endowed with several metrics. 
In this article, we concentrate on one of them, called bottleneck distance. 
\begin{defn}\label{def_bottdist}
Let $\p$ and $\psi$ be continuous functions. 
The \textbf{bottleneck distance} between the persistence diagrams of $\p$ and $\psi$ is
\[
d_\textup{B}(\mathrm{Dgm}(\p), \mathrm{Dgm}(\psi))=\inf _{\sigma} \mathrm{cost}(\sigma),
\]
where $\sigma$ varies over all the multiset bijections, also called \textbf{matchings}, from $\mathrm{Dgm}(\p)$ to $\mathrm{Dgm}(\psi)$, and $\mathrm{cost}(\sigma)=\sup_{p\in \mathrm{Dgm}(\p)} d(p, \sigma(p))$.
\end{defn}

In \cite{Cohen_Steiner_stability}, a stability theorem for the bottleneck distance is proven:
\begin{thm}
\label{thm:bottleneck_stability}
    Let $X$ be a triangulable space and $\p, \s \colon X \to \R$ two continuous functions satisfying Assumption~\ref{jdshdsds}. Then
    $\dmatch{\textup{Dgm}(\p)}{\textup{Dgm}(\s)} \le \norm{\p-\s}{\infty}$.
\end{thm}

As a consequence of Theorem \ref{thm:bottleneck_stability}, since our functions $\p, \psi\colon X\to \R$ are assumed to be continuous with $X$ compact, the distance $\lVert\p-\psi\rVert_\infty$ must be finite. Thus, $\dmatch{\textup{Dgm}(\p)}{\textup{Dgm}(\s)}$ must also be finite.

If there exists $\overline{\sigma}$ such that 
$d_\textup{B}(\mathrm{Dgm}(\p), \mathrm{Dgm}(\psi))= \textup{cost}(\overline{\sigma})$
, we refer to it as an \textbf{optimal matching}. 
By Assumption~\ref{jdshdsds} and the compactness of (the supports of) persistence diagrams, an optimal matching always exists
(see Propositions \ref{cpt} and \ref{opt_matching} in Appendix A).
Furthermore, one can show that the bottleneck distance is an extended metric between persistence diagrams, or equivalently between PBNFs, and an extended pseudo-metric between functions.

\subsection{$2$-parameter Persistent Homology}\label{subsec:biparameter}
Persistent homology can be also defined for filtering functions of the form $\p\colon X\to \R^n$, with $n\ge 2$. 
In this case, we can still define PBNFs, but not persistence diagrams, as in the $1$-parameter case. 
In what follows, we concentrate on the case when $n=2$. 
This is special because it is possible to reduce the information condensed in the PBNF of the function to, not one, but a collection of persistence diagrams. 

Let $X$ be a compact topological space and $\varphi \colon X \to \R^2$, with $\p(x)=(\p_1(x), \p_2(x))$, a continuous function.
Analogously to the $1$-parameter case, every $(u_1,u_2)$ in $\mathbb{R}^2$ determines a sublevel set $(\varphi\le(u_1,u_2)) = \{x \in X \mid \varphi_1(x) \leq u_1, \; \varphi_2(x) \leq u_2\}$.
The inclusion $(\varphi\le(u_1,u_2))\subset (\varphi\le(v_1,v_2))$, for $(u_1,u_2)\le (v_1,v_2)$, induces a linear map in homology $\iota^*_{(u_1,u_2),(v_1,v_2)} \colon H_k(\varphi\le(u_1,u_2)) \to H_k(\varphi\le(v_1,v_2))$, for every integer $k$. 
As before, we will assume that:
\begin{ass}\label{fhwaoijdncks}
$H_k(\varphi\le(u_1,u_2))$ is a finitely generated $\mathbb{K}$-vector space for every $k\in \mathbb{Z}$ and for every $(u_1,u_2)\in \R^2$. 
\end{ass}
The $k$\textbf{-th persistent homology group} of $\p$ at $((u_1, u_2),(v_1, v_2))$ is defined as the subspace $H^\p_k((u_1,u_2),(v_1,v_2)) = \textup{Im }\iota^*_{(u_1,u_2),(v_1,v_2)}\subset H_k(\p \le (v_1,v_2))$. 
The \textbf{$k$-persistent Betti numbers function}, or PBNF, of $\p$ at $((u_1,u_2),(v_1,v_2))$, $\beta^\p_k((u_1,u_2),(v_1,v_2))$, is defined as $\textup{dim } H^\p_k((u_1,u_2),(v_1,v_2))$.

\begin{figure}
    \centering
\includegraphics[width=6cm]{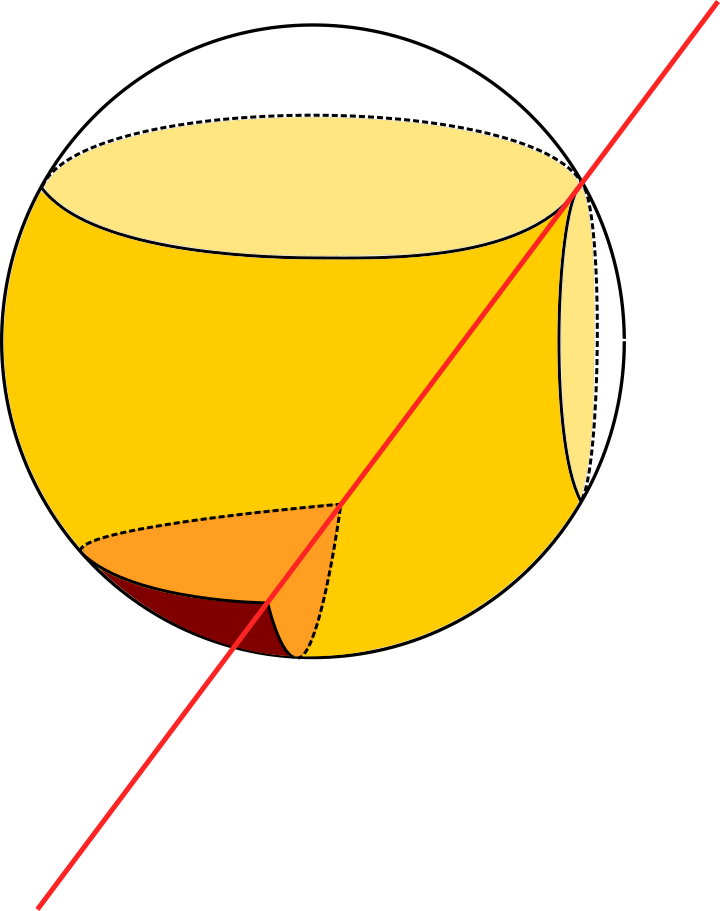}
    \caption{Two different sublevel sets of the projection of a sphere in $\R^3$ onto the plane $y=0$. 
    Note that the smaller sublevel set is contractible, whereas the bigger sublevel set has the homotopy type of $ S^1$.}
    \label{fig:filt}
\end{figure}

A convenient way to study $2$-parameter persistence is by reducing it to the $1$-parameter case. 
To do this, we consider the set of all lines, $r_{(a,b)}$, in $\R^2$ with positive slope parameterised by $]0,1[\times \R$ as $r_{(a,b)} = \{t(a,1-a)+(b,-b) \mid t \in \R\}$.
Each line $r_{(a,b)}$ induces a $1$-parameter filtration associating each point $(u_1(t),u_2(t))$ in $r_{(a,b)}$ with the sublevel set $(\p\le (u_1(t),u_2(t)))$, for every $t$ in $\R$.
The subspaces in this filtration coincide with those that arise from the function $\p_{(a,b)}\colon X\to \R$ defined as $\p_{(a,b)}(x) = \max \{\frac{\p_1(x)-b}{a}, \frac{\p_2(x)+b}{1-a}\}$.
Note that, since $H_k(\varphi\le (u_1, u_2))$ satisfies Assumption~\ref{fhwaoijdncks} for every $k$ and for every $(u_1,u_2)$, then the restriction $\p_{(a,b)}$ satisfies Assumption~\ref{jdshdsds}.
We will also
consider the normalised function $\p^*_{(a,b)}(x) = \min \{a,1-a\} \p_{(a,b)}(x)$, since it
will be of use in what follows.
For details about this method, known as \emph{foliation method}, we refer to \cite{pers_Betti_stable}.
Figure~\ref{fig:filt} shows an example of a sublevel set filtration along a line of positive slope of a sphere in $\R^3$. 
Taking advantage of this method to represent $2$-parameter filtrations, it is possible to define the following metric between PBNFs of $\R^2$-valued filtering functions. 

\begin{defn}
Consider two continuous functions $\p,\psi\colon X\to \R^2$, the \textbf{matching distance} between them is defined as  
\[
\Dmatch{\p}{\psi} = \sup_{(a,b) \in ]0,1[\times\R} \dmatch{\dgm{\p}{a}{b}}{\dgm{\psi}{a}{b}}.
\]
If there exists $(\overline a, \overline b) \in ]0,1[ \times \R$ such that $\Dmatch{\p}{\psi} =\dmatch{\dgm{\p}{\overline a}{\overline b}}{\dgm{\psi}{\overline a}{\overline b}}$ we say that $(\overline a, \overline b)$ \textbf{realises} $\Dmatch{\p}{\psi}$.
\end{defn}

Analogously to the bottleneck distance, the matching distance can be seen as an extended metric between PBNFs and an extended pseudo-metric between filtering functions. 
The notation used here recalls only the functions, however, we mean it in both senses.

The matching distance is stable, in the following sense (see \cite[Theorem 4.4]{pers_Betti_stable}):
\begin{thm}\label{stability}
Let $X$ be a triangulable space and $\p,\psi\colon X\to \R^2$ two continuous functions satisfying Assumption~\ref{fhwaoijdncks}. Then
$\Dmatch{\p}{\psi} \leq \norm{\p - \psi}{\infty}$.
\end{thm}

Theorem~\ref{stability} implies that $\Dmatch{\p}{\s}$ is finite under our assumptions that $\p$ and $\psi$ are continuous functions on a compact domain.

As a pseudo-metric between filtering functions, the matching distance can be interpreted as an approximation of the $L^{\infty}$-distance between the same functions. 
We now present two results showing a relation between the parameter values in $]0,1[\times\R$ where the matching distance can be realised and how well the matching distance approximates the $L^{\infty}$-distance. 
This partially motivates our interest in a qualitative analysis of the parameter space of the filtering lines involved in the definition of the matching distance and, particularly, our attention to lines of slope 1 or, equivalently, with $a=\frac{1}{2}$ in our parameterisation.  

\begin{thm}\label{snfvjhfhiksld}
Let $\p =(\p_1,\p_2), \psi=(\psi_1, \psi_2) \colon X \to \R^2$ be continuous functions.
    Let $(\overline a, \overline b)$ realise $\Dmatch{\p}{\s}$.
    If $\norm{\p - \s}{\infty} - \Dmatch{\p}{\s} < \norm{\p_2 - \s_2}{\infty} - \norm{\p_1 - \s_1}{\infty}$, then $\overline a  > \frac{\norm{\p_1 - \s_1}{\infty}}{\norm{\p_1 - \s_1}{\infty}+\norm{\p_2 - \s_2}{\infty}}$.
\end{thm}
\begin{proof}
Denote by $c$ the difference $ \norm{\p - \s}{\infty} - \Dmatch{\p}{\s}$, which is greater or equal to zero by Theorem~\ref{stability}.
    Observe that $\frac{\norm{\p_1 - \s_1}{\infty}}{\norm{\p_1 - \s_1}{\infty}+\norm{\p_2 - \s_2}{\infty}} < \frac{1}{2}$,
    and therefore the result holds for $\overline a \ge \frac{1}{2}$.
Thus, let $\overline a < \frac{1}{2}$.
    Applying the inequality $\lvert \max\{A, B\}- \max\{C,D\}\rvert \le \max\{\lvert A- C\rvert, \lvert B-D\rvert\}$, where $A,B,C$ and $D$ are real numbers, we obtain
    \begin{align*}
        \norm{\p^*_{(\overline a, \overline b)} - \s^*_{(\overline a, \overline b)}}{\infty} &\le \min\{\overline a,1-\overline a\} \max  \bigg\{\frac{\norm{\p_1 - \s_1}{\infty}}{\overline a}, \frac{\norm{\p_2 - \s_2}{\infty}}{1-\overline a}\bigg\} \\
        &= \overline a \max  \bigg\{\frac{\norm{\p_1 - \s_1}{\infty}}{\overline a}, \frac{\norm{\p_2 - \s_2}{\infty}}{1-\overline a}\bigg\} \\
        &= \max  \bigg\{\norm{\p_1 - \s_1}{\infty}, \frac{\overline a}{1-\overline a}\norm{\p_2 - \s_2}{\infty}\bigg\}.
    \end{align*}
Suppose by contradiction that $\overline a \le \frac{\norm{\p_1 - \s_1}{\infty}}{\norm{\p_1 - \s_1}{\infty}+\norm{\p_2 - \s_2}{\infty}}$.
Then $\norm{\p_1 - \s_1}{\infty} \ge \frac{\overline a}{1-\overline a}\norm{\p_2 - \s_2}{\infty}$ and, hence, $\norm{\p^*_{(\overline a, \overline b)} - \s^*_{(\overline a, \overline b)}}{\infty} \le \norm{\p_1 - \s_1}{\infty}$.
As a consequence, $\norm{\p_1 - \s_1}{\infty}\ge\Dmatch{\p}{\s} = \norm{\p - \s}{\infty} - c=\norm{\p_2 - \s_2}{\infty} - c$, contradicting the hypothesis $c < \norm{\p_2 - \s_2}{\infty} - \norm{\p_1 - \s_1}{\infty}$. 
    Therefore, $\overline a  > \frac{\norm{\p_1 - \s_1}{\infty}}{\norm{\p_1 - \s_1}{\infty}+\norm{\p_2 - \s_2}{\infty}}$.
\end{proof}
\begin{cor}
    Let $\p, \psi \colon X \to \R^2$ be continuous functions.
    Let $(\overline a, \overline b)$ realise $\Dmatch{\p}{\s}$. Assume $\Dmatch{\p}{\s} = \norm{\p - \s}{\infty}$.
    If the inequality $\norm{\p_2 - \s_2}{\infty} > \norm{\p_1 - \s_1}{\infty}$ holds, then $\overline a  > \frac{\norm{\p_1 - \s_1}{\infty}}{\norm{\p_1 - \s_1}{\infty}+\norm{\p_2 - \s_2}{\infty}}$.
\end{cor}

\subsection{Extended Pareto grid}
In this section we recall the relation between a differential construction associated with a smooth $\R^2$-valued function, called the extended Pareto grid, and the points of the persistence diagrams $\dgm{\p}{a}{b}$. This connection is established in the Position Theorem proven in~\cite{coher_match}.

Let $M$ be a closed smooth Riemannian manifold and $\p=(\p_1,\p_2)\colon M\to \mathbb{R}^2$ a smooth function. 

\begin{defn}\label{def:prop_functions}
The \textbf{Jacobi set} of $\p$ is defined by 
$\mathbb{J}(\p)=\{p\in M\mid \nabla \p_1=\lambda \nabla \p_2 \text{ or }\nabla \p_2=\lambda \nabla \p_1, \text{ for some }\lambda\in \mathbb{R}\}$.
The \textbf{Pareto set} of $\p$ is the subset of $\mathbb{J}(\p)$ given by
$\mathbb{J}_P(\p)=\left\{p\in M \mid \nabla \p_1=\lambda \nabla \p_2 \text{ or }\nabla \p_2=\lambda \nabla \p_1, \text{ for some }\lambda\le 0\right\}$.
\end{defn}

Unless differently specified, in the rest of the article, the smooth filtering functions considered will satisfy the following conditions.
\begin{ass}\label{ewaiowjnvornsje}
The smooth function $\p=(\p_1, \p_2)\colon M\to \R^2$ satisfies the following conditions.
\begin{enumerate}[label=(\roman*)]
    \item No point $p$ exists in $M$ at which both $\nabla \p_1$ and $\nabla \p_2$ vanish.
    \item $\mathbb{J}(\p)$ is a $1$-manifold smoothly embedded in $M$ consisting of finitely many components, each one diffeomorphic to a circle.
    \item $\mathbb{J}_P(\p)$ is a $1$-dimensional closed submanifold of $\mathbb{J}(\p)$ with boundary.
    \item If we denote by $\mathbb{J}_C(\p)$ the subset of $\mathbb{J}(\p)$ where $\nabla \p_1$ and $\nabla \p_2$ are orthogonal to $\mathbb{J}(\p)$, then the connected components of $\mathbb{J}_P(\p)\setminus\mathbb{J}_C(\p)$ are finite in number, each one being diffeomorphic to an interval.
    With respect to any parameterisation of each component, one of $\p_1$ and $\p_2$ is strictly increasing and the other is strictly decreasing. Each component can meet critical points for $\p_1$, $\p_2$ only at its endpoints.
\end{enumerate}
\end{ass}

 These properties are generic in the set of smooth functions from $M$ to $\R^2$ (see \cite{wan}).

In the following, we will consider $\R \cup \{\infty\}$, denoted by $\Rinf$, with the topology generated by $]s,t[$, with $s < t$ in $ \R$, and $]s,\infty]$.

\begin{defn}
\label{def:epg}
Denote by $\{p_1,\dots, p_h\}$ and $\{q_1,\dots, q_k\}$, respectively, the critical points of $\p_1$ and $\p_2$. 
The \textbf{extended Pareto grid} of $\p$ is defined as the union 
\[
\Gamma (\p)=\p\left(\mathbb{J}_P(\p)\right)\cup \left(\bigcup _i v_i\right) \cup  \left(\bigcup _j h_j\right)
\]
where $v_i$ is the vertical half-line $\{(x,y)\in \Rinf^2\mid x=\p_1(p_i), y\ge \p_2(p_i)\}$ and
$h_j$ is the horizontal half-line $\{(x,y)\in \Rinf^2\mid x\ge \p_1(q_j), y= \p_2(q_j)\}$.
We refer to these half-lines as \textbf{improper contours} and to the closure of the images under $\p$ of the connected components of $\mathbb{J}_P(\p)\setminus\mathbb{J}_C(\p)$ as \textbf{proper contours}. 
The set of all contours of $\Gamma(\p)$ is denoted by $\mathrm{Ctr}(\p)$.
\end{defn}

\begin{figure}
\includegraphics[width=1\textwidth]{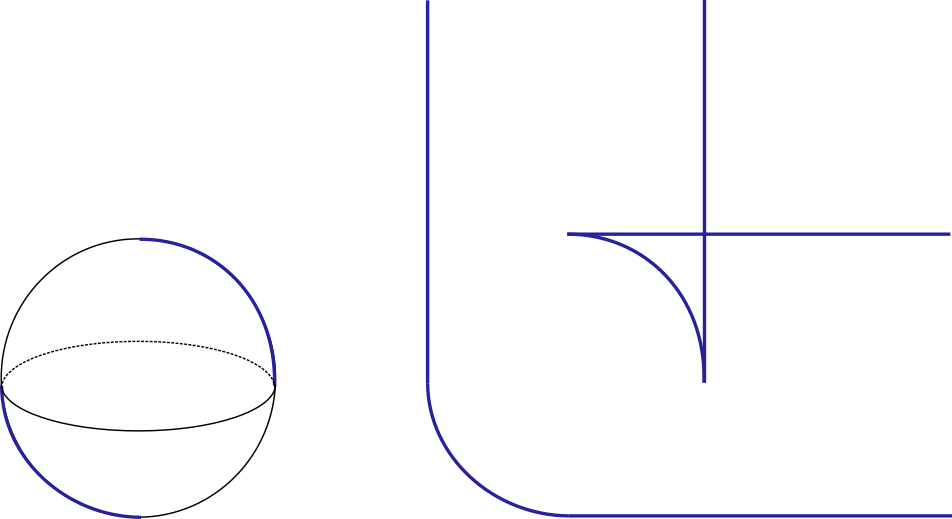}
    \caption{Consider the sphere $S^2$ in $\R^3$ and its canonical projection $\pi$ on the plane $y=0$. 
    On the left, the subset in blue is the Pareto set of $\pi$. 
    On the right, the extended Pareto grid of $\pi$.
    }
    \label{fig:EPG_sphere}
\end{figure}

For an example of extended Pareto grid of a function on a surface, see Figure~\ref{fig:EPG_sphere}.

\begin{rem}\label{hijjdfshijs} The following are consequences, respectively, of (ii), (iii), (iv) and, again, (iv) in Assumption \ref{ewaiowjnvornsje}.
\begin{itemize}
\item Ctr$(\p)$ contains a finite number of contours. 
\item The closure in $\mathbb J (\p)$  of every proper contour is diffeomorphic to a closed interval of $\R$.
\item Every contour can be parameterised as a curve whose two coordinates are respectively non-increasing and non-decreasing.
In particular, this implies that the tangent line to any proper contour at any point of its interior cannot have strictly positive slope.
\item For every $(a,b)\in ]0,1[\times \R$ and for every proper contour $\alpha\in \mathrm{Ctr}(\p)$, the intersection of $r_{(a,b)}$ and $\alpha$ is either empty or a single point.
\end{itemize}
\end{rem}

\begin{rem}
Note that the Definition \ref{def:epg} is an extension of the original definition given in \cite{coher_match}.
In the current definition the improper contours are diffeomorphic to intervals. 
This will allow us to introduce a suitable notion of extended intersection in Section \ref{sec_inters}.
\end{rem}

\subsection{Preliminary results}
\label{sec:neededresults}

In~\cite{slope1} it was observed that persistence diagrams can also be associated with $\p_{(a,b)}^*$ when $a=0, 1$. 
In particular, as a consequence of~\cite[Theorem 4.3]{slope1}, the function $(a,b)\mapsto\p^*_{(a,b)}$ defined on $]0,1[\times\R$ can be extended to $[0,1]\times \R$ as 
\begin{align*}
\p^*_{(0, b)}(x) & =\lim_{(a',b') \to (0,b)} \p^*_{(a',b')}=\max\{\p_1(x)-b,0\},\\
\p^*_{(1, b)}(x) & =\lim_{(a',b') \to (1,b)} \p^*_{(a',b')}=\max\{0,\p_2(x)+b\}.
\end{align*}
Thus, given two functions $\p$ and $\psi$, the persistence diagrams $\mathrm{Dgm}(\p_{(a,b)}^*)$ and $\mathrm{Dgm}(\psi_{(a,b)}^*)$ are defined for every $(a,b)$ in $[0,1]\times \R$, and the function
$$
(a,b)\mapsto \dmatch{\dgm{\p}{a}{b}}{\dgm{\s}{a}{b}}
$$
can also be extended to $[0,1]\times \R$. 
Moreover, such extension is continuous. 
In particular, 
\begin{equation}\label{lemma_stability_position}
\lim_{a \to 0} \dmatch{\dgm{\p}{a}{b}}{\dgm{\p}{0}{b}}=0, \;\;\;\;
\lim_{a \to 1} \dmatch{\dgm{\p}{a}{b}}{\dgm{\p}{1}{b}}=0.
\end{equation}

In \cite[Proposition 4.4]{slope1} it was also showed that the matching distance can be realised by parameter values lying in a bounded region of $[0,1]\times\R$: 

\begin{prop}\label{slope1_Dmatch_realises}
    Let $\overline C = \max \{\norm{\p}{\infty},\norm{\s}{\infty}\}$.
    There exists $(\overline a, \overline b)\in[0,1]\times[-\overline C, \overline C]$ such that 
$\Dmatch{\p}{\s} = \dmatch{\dgm{\p}{a}{b}}{\dgm{\s}{a}{b}}$.
\end{prop}

The previous result allows us restrict to the compact space of parameters, $[0,1]\times[-\overline C,\overline C]$.

\section{Intersection operator and Position Theorem}\label{sec_inters}

In this section we present a definition of intersection between filtering lines and contours of the extended Pareto grid. In particular, we introduce a notion of intersection in $\R^2$ that does not coincide with the usual one. This definition is paramount to extend the Position Theorem of \cite{coher_match} to the case of vertical and horizontal filtering lines.

Before proceeding we need to define the filtering line $r_{(a,b)}$ for $a=0,1$. In particular, we have that $r_{(0,b)}$ is the horizontal line $ \{x = b\}$ and $r_{(1,b)}$ is the vertical line $\{y= - b\}$.
Since the functions we are considering satisfy property (iv) in Assumption~\ref{ewaiowjnvornsje} and the filtering lines $r_{(a,b)}$, with $(a,b)\in [0,1]\times \R$, have non-negative slope, the canonical intersection between a contour and a filtering line can be empty, a single point or infinitely many points (see Remark~\ref{hijjdfshijs}). 

\begin{defn}\label{def:intersection_operator}
    For every $(a,b) \in [0,1]\times \R$ and every $\alpha \in \mathrm{Ctr}(\p)$ set
\[
    \Rinf^2 \supseteq r_{(a,b)} \hcap \alpha =
    \begin{cases}
        r_{(a,b)} \cap \alpha \quad &\textup{if } a \ne 0,1, \\
        \{ P = \lim_{n \to \infty}P_n \mid P_n \in r_{(a_n,b_n)} \cap \alpha\} &\textup{if } a \in \{0,1\}.
    \end{cases}
\]
where $((a_n,b_n))_{n \in \N}$ is a sequence
in $]0,1[\times\R$
such that $\lim_{n\to \infty} a_n = a$ and $\lim_{n\to \infty} b_n = b$.

    The collection of points in $r_{(a,b)} \hcap \alpha$, when $\alpha$ varies in $\mathrm{Ctr}(\p)$, is denoted by $r_{(a,b)} \hcap \Gamma(\p)$. 
\end{defn}

Note that if $r_{(a,b)} \cap \alpha$ is unique then $r_{(a,b)} \hcap \alpha$ is unique. 
In this case, the new intersection operator, $\hcap$, coincides with the usual one. 
Moreover, $r_{(a,b)} \hcap \alpha$ can coincide with $r_{(a,b)} \cap \alpha$ only when $\alpha \subset r_{(a,b)}$.
The new intersection and the canonical one are different when the filtering line is parallel to, but does not contain, the improper contour. 
For example, if $\alpha$ is a vertical improper contour and $a=0$, $r_{(a,b)} \hcap \alpha$ is empty if $\alpha$ is on the left of $r_{(a,b)}$, and it is only one point if $\alpha$ is on the right of $r_{(a,b)}$.

\begin{ex}

    Fix $b\in[-\overline C,\overline C]$.
    Recall that the line $r_{(a,b)}$ is given by the parametric equation $(x(t),y(t)) = (at+b,(1-a)t-b)$.
    Let $h_{x_0,y_0}$ be the improper horizontal contour given by $\{x \ge x_0, \; y = y_0\}$, with $x_0 >b, y_0 > -b$. 
    There exists a strictly positive real value $A<1$ such that for every $a \in [A,1[$ the filtering line $r_{(a,b)}$ intersects the contour, and for every $a\in ]0,A[$ it does not. 
    In particular, $A$ is uniquely determined by the condition $(x_0,y_0) \in r_{(A,b)}$, that holds if, and only if, $A = \frac{x_0-b}{x_0+y_0}$.
    Note that $x_0+y_0 > 0$ because $x_0 >b , y_0>-b$.
Consider the sequence $(a_n,b)=\left(\left(\frac{n-1+A}{n},b\right)\right)_{n\in\N}\subset [A,1[ \times \{b\}$ converging to $(1,b)$.
     It yields the sequence $(P_n)_{n \in \N}$ of points in $r_{(a_n,b)}\cap h_{x_0,y_0}$
     with ordinate $y_0=\frac{1-A}{n}t - b$. 
     Therefore, $x_n = \left(\frac{n-1+A}{n}\right)t+b=(n-1+A)\frac{b+y_0}{1-A} + b$, and $\lim_{n \to \infty} x_n = \infty$, allowing us to conclude that $r_{(1,b)} \hcap h_{x_0,y_0} = \{(\infty, y_0)\}$.
\end{ex}

\begin{figure}[htp]
    \hspace*{\fill}%
    \raisebox{-\height}{\includegraphics[width=0.4\textwidth]{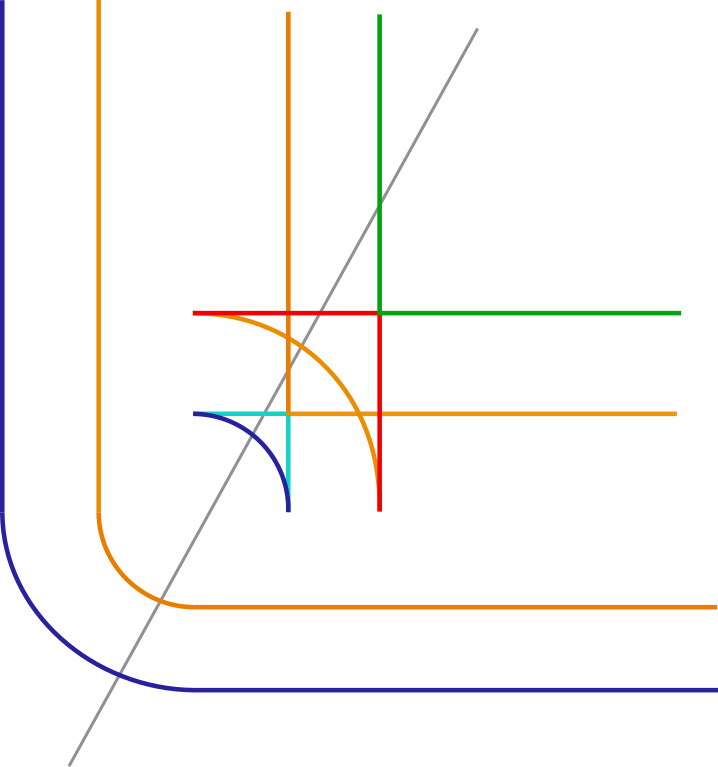}}%
    \hfill
    \raisebox{-\height}{\includegraphics[width=0.4\textwidth]{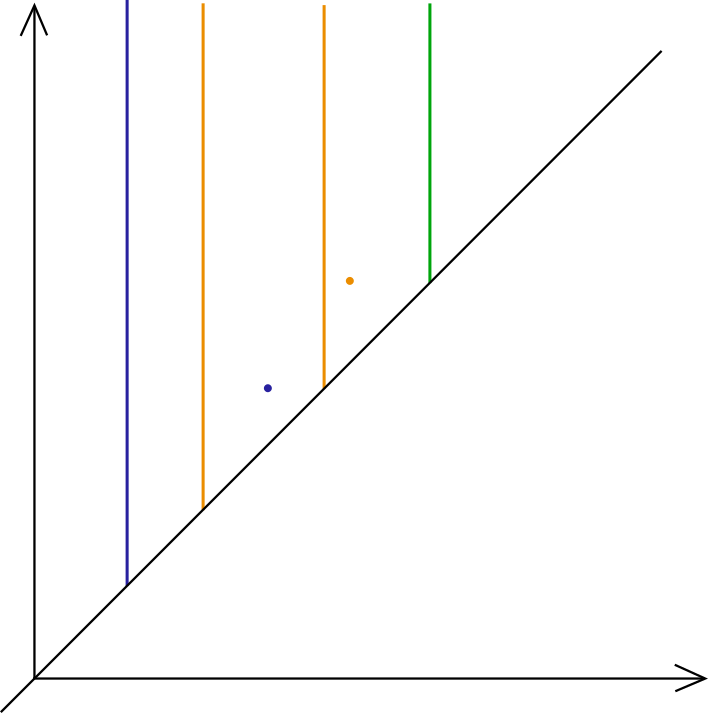}}%
    \hspace*{\fill}
    \caption{Consider a torus in $\R^3$ and a suitable projection $\pi$ on $\R^2$. 
    On the left, the extended Pareto grid of $\pi$, together with a filtering line $r_{(a,b)}$ in grey. 
    The colour choice for the contours is as follows: blue/cyan are for the birth/death of a connected component, orange/red are for the birth/death of a 1-cycle and green is for the birth of a 2-cycle. 
    The figure on the right is a superposition of the persistence diagrams of $\pi^*_{(a,b)}$ in degree 0,1 and 2, in blue, orange and green, respectively. There the half-lines represent the essential cornerpoints}
    \label{fig:position_thm_torus}
\end{figure}

In~\cite[Theorem 2]{coher_match}, the Position Theorem was stated and proven. 
It gives a correspondence of the points of the persistence diagram of $\p_{(a,b)}^*$ with the intersections of the line $r_{(a,b)}$ with the extended Pareto grid of $\p$, for every $(a,b)$ in $]0,1[\times \R$. 
Figure~\ref{fig:position_thm_torus} shows an example of this for a projection of a torus in $\R^3$.
We show a generalisation of this theorem that comprehends parameter values in $[0,1]\times \R$, providing the finite coordinates of the points in the persistence diagrams associated with $\p_{(0,b)}^*$ and $\p_{(1,b)}^*$ for every $b$ in $\R$. 

\begin{thm}\label{position_generalised}
    Let $(a,b) \in [0,1] \times \R$ and $p \in \dgm{\p}{a}{b} \setminus \{\Delta\}$.
    For each finite coordinate $w$ of $p$, a point $P = (x, y) \in r_{(a,b)} \hcap \Gamma(\p)$ exists such that 
    
    $$
    w = \begin{cases}
        \min\left\{1,\frac{1-a}{a}\right\}(x-b) \quad\textup{ if } a \in [0,1[, \\
        \min\left\{1,\frac{a}{1-a}\right\}(y+b) \quad\textup{ if } a \in ]0,1].
    \end{cases}
    $$
    
    with the conventions $\frac{1}{0} = \infty$, $\min\{1, \infty\} = 1$.
\end{thm}
\begin{proof}
For $a \in ]0,1[$, the proof is the original Position Theorem as it appears in \cite{coher_match}.
Let $a = 0$ and $b \in \R$. Let $p \in \dgm{\p}{0}{b}$ and $w$ be a finite coordinate of $p$.
The first equality in Equation (\ref{lemma_stability_position}) guarantees that we can choose a sequence $((a_n, b))_{n \in \N}$ 
in $]0,1[ \times \R$ converging to $(0,b)$ in $[0,1] \times \R$ and, by possibly extracting a subsequence, we can assume that there exists a sequence of cornerpoints
\[
    p_n \in \dgm{\p}{a_{n}}{b}, \textup{ with } d(p_n,p) < \frac{1}{n}.
\]

Furthermore, we can assume that $a_{n} < \frac{1}{2}$ for any index $n$, again up to subsequences. 
Let $w_n$ be the finite coordinate of $p_n \in \dgm{\p}{a_{n}}{b}$ converging to the coordinate $w$ of $p \in \dgm{\p}{0}{b}$. 
Applying the classical Position Theorem to $w_n$, we have that there exists a point $P_n\in r_{(a_n,b)}\cap \Gamma(\varphi)$ such that
\[
w_n = \min \left\{1, \frac{1-a_n}{a_n}\right\} (x_{P_n} - b) = x_{P_n} - b.
\]

By compactness of $\Gamma(\p)$ in $\Rinf^2$, the sequence $(P_n)_{n \in \N}$, up to subsequences, converges to a point $P \in \Gamma(\p)$.
Then, we have that
\[w = \lim_{n \to \infty} w_n = \lim_{n \to \infty} x_{P_{n}} - b = x_P - b,\] 
as we wanted to show.

The proof for the case $a = 1$ is analogous.
\end{proof}

\begin{rem}\label{uwuhhddijdsfjdsfhi}
Recall that the bottleneck distance depends on the definition of the metric $d$ (see (\ref{opcklfajkclòkweaa})).
For any $(a,b)\in [0,1]\times[-\overline{C}, \overline{C}]$, if $\dmatch{\dgm{\p}{a}{b}}{\dgm{\s}{a}{b}}>0$ then it is either equal to $\lvert w-w'\rvert$ or to $\frac{1}{2}\lvert w-w'\rvert$, where $w$ and $w'$ are finite coordinates of points in $\dgm{\p}{a}{b}\cup \dgm{\s}{a}{b}$.
In particular, the coefficient $\frac{1}{2}$ is obtained when matching a cornerpoint to the diagonal $\Delta$. 
Theorem~\ref{position_generalised} ensures that there exist $P,Q\in (\Gamma(\p)\cup \Gamma(\s))\hat \cap r_{(a,b)}$ such that  
    \[  
    w=\min\left\{1,\frac{1-a}{a}\right\}(x_P-b) = \min\left\{1,\frac{a}{1-a}\right\}(y_P+b),
    \]
    \[
    w'=\min\left\{1,\frac{1-a}{a}\right\}(x_Q-b) = \min\left\{1,\frac{a}{1-a}\right\}(y_Q+b).
    \]
    Thus, if $a\le \frac{1}{2}$, the bottleneck distance is either $\lvert x_P-x_Q\rvert$ or $\frac{1}{2}\lvert x_P-x_Q\rvert$, and  if $a\ge \frac{1}{2}$, the bottleneck distance is either $\lvert y_P-y_Q\rvert$ or $\frac{1}{2}\lvert y_P-y_Q\rvert$. 
    In both cases, we say that $P$ and $Q$ in $(\Gamma(\p)\cup \Gamma(\s))\cap r_{(a,b)}$ \textbf{realise} the bottleneck distance between $\dgm{\p}{a}{b}$ and $\dgm{\s}{a}{b}$. 

\end{rem}

\section{Special and ultraspecial values}
\label{sec:special}

Consider the filtering functions $\p$ and $\psi$ and their respective extended Pareto grids, $\Gamma(\p)$ and $ \Gamma(\psi)$. 
There are pairs $(a,b) \in ]0,1[ \times [-\overline C,\overline C]$ in which the optimal matching between $\dgm{\p}{a}{b}$ and $\dgm{\s}{a}{b}$ may change abruptly, as the following example shows.

\begin{ex}\label{ex_special}
Consider the improper contours of $\p$ and $\s$ (in blue and red) in Figure \ref{fig:tikz_ex_special}. 
Let $a \le \frac{1}{2}$.
The three parallel filtering lines, in grey, $r_{(a,b')}$, $r_{(a,b)}$ and $r_{(a,b'')}$, intersect the contours, hence, by Theorem \ref{position_generalised}, they give the $x$-coordinates of some cornerpoints $A,B,C, C',C''$.
The persistence diagrams associated with the different filtering lines are also represented in Figure \ref{fig:tikz_ex_special}.
Since the $x$-coordinates of $A,B$ are given by a vertical contour and $a \le \frac{1}{2}$, these two points remain fixed for each $u \in \{b,b',b''\}$, up to a translation parallel to the diagonal.
Instead, the intersections with the horizontal improper contour with $r_{(a,b')}$, $r_{(a,b)}$ and $r_{(a,b'')}$ give the translation of the $x$-coordinate of $C'$, $C$ and $C''$.
Note that for simplicity we are assuming that the cost of the optimal matching does not depend on any other cornerpoint.
Assuming that these cornerpoints are sufficiently far from the diagonal, when the filtering line is translated from $r_{(a,b')}$ to $r_{(a,b'')}$, the optimal matching changes abruptly at $(a,b)$, where there are two possible optimal matchings realising the bottleneck distance, one sending $A$ to $C$ and $B$ to the diagonal, and another one sending $B$ to $C$ and $A$ to the diagonal. 

\begin{figure}
        \includegraphics[width=\textwidth]{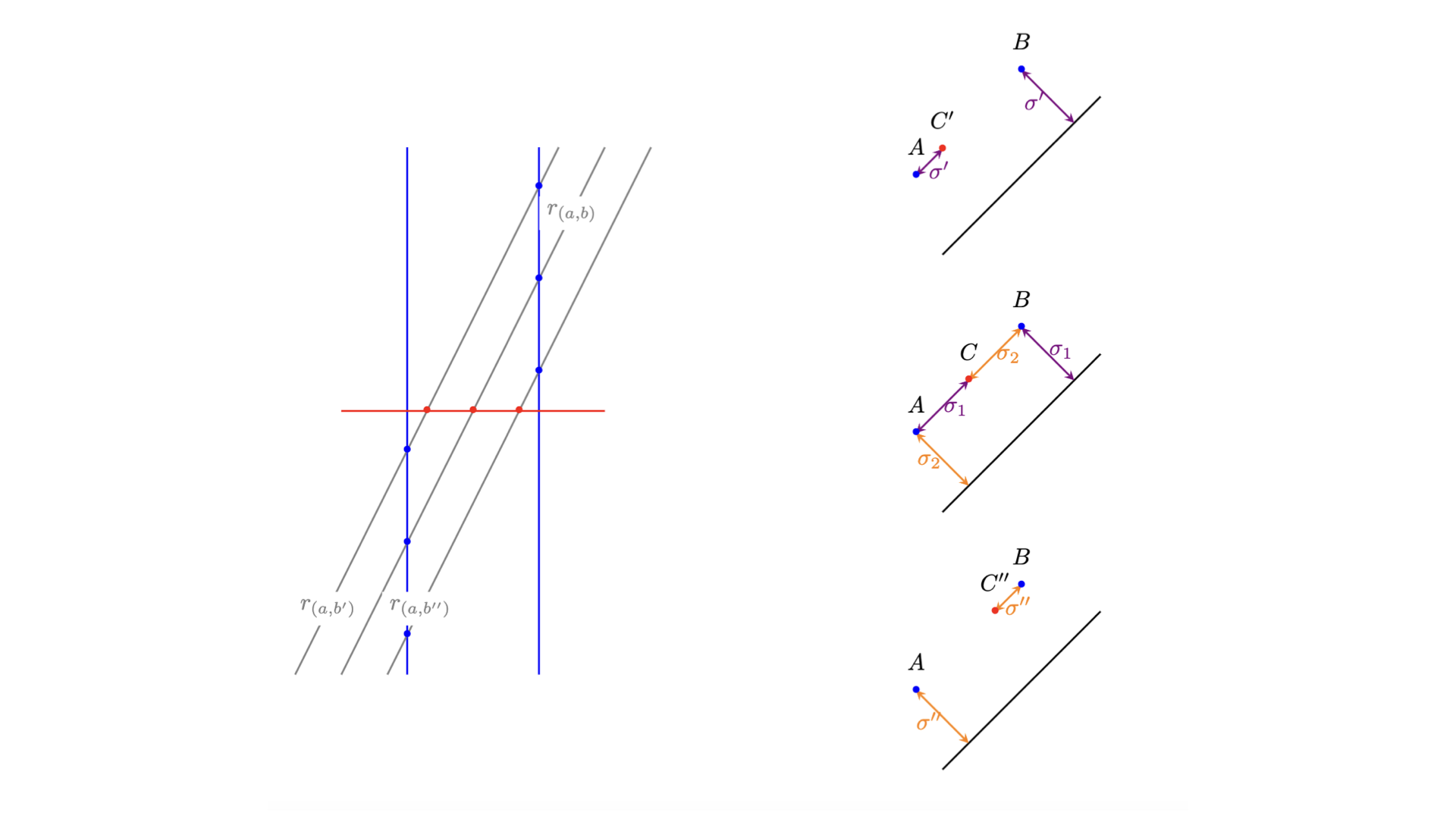}
        \caption{On the left, $\Gamma(\p)$ in blue, $\Gamma(\s)$ in red and the filtering lines $r_{(a,b)}$, $r_{(a,b')}$ and $r_{(a,b'')}$. On the right, the persistence diagrams associated with these filtering lines along with their corresponding optimal matchings. 
        Note that for $r_{(a,b)}$ there are two possible optimal matchings.}
        \label{fig:tikz_ex_special}
    \end{figure}
\end{ex}

The phenomenon described in Example~\ref{ex_special} motivates the definition of special set introduced in \cite{slope1}. 

\begin{defn}\label{spec}
Let $\textup{Ctr}(\p, \psi) = \textup{Ctr}(\p) \cup \textup{Ctr}(\psi)$ be the set of all curves that are contours of $\p$ or $\psi$. 
The \textbf{special set}
of $(\p,\psi)$, denoted by $\spe{\p}{\psi}$, is the
collection of all $(a,b)$ in $]0,1[ \times [-\overline{C},\overline{C}]$ for which two different pairs $\{\alpha_P, \alpha_Q\}$,
$\{\alpha_R, \alpha_S\}$ of contours in $\mathrm{Ctr}(\p, \psi)$ intersecting $r_{(a,b)}$
exist, such that 
    \begin{itemize}
        \item $c_1\lvert x_P-x_Q\rvert=c_2\lvert x_{R}-x_{S}\rvert$, with $c_1,c_2\in \{1,2\}$, if $a\le \frac{1}{2}$,
        \item $c_1\lvert y_P-y_Q\rvert=c_2\lvert y_{R}-y_{S}\rvert$, with $c_1,c_2\in \{1,2\}$, if $a\ge \frac{1}{2}$,
    \end{itemize}
    where $P=P_{(a,b)}=r_{(a,b)}\cap \alpha_P$,
    $Q=Q_{(a,b)}=r_{(a,b)}\cap \alpha_Q$, $R=R_{(a,b)}=r_{(a,b)}\cap \alpha_R$ and
    $S=S_{(a,b)}=r_{(a,b)}\cap \alpha_S$, and $x_{*}$, $y_{*}$ denote abscissas and ordinates of these points, respectively.
    An element of the special set $\spe{\p}{\psi}$ is called a \textbf{special value} of the pair $(\p, \psi)$.
    Any two pairs of contours for which one of the conditions above is satisfied is called $(a,b)$\textbf{-special}.
    
\end{defn}
\begin{rem}
    Note that, in the above definition, 
    the contours $\alpha_P, \alpha_Q,\alpha_R, \alpha_S$ do not necessarily need to differ.
    For example, the definition allows one among $\alpha_P = \alpha_R$, $\alpha_P = \alpha_S$, $\alpha_Q=\alpha_R$, $\alpha_Q = \alpha_S$ to hold, but not two of these conditions simultaneously.
    Moreover, all these contours may belong to the same extended Pareto grid.
\end{rem}

\begin{rem}\label{special_double_point}
    Let $r_{(a,b)}$, with $(a, b) \in ]0,1[\times [-\overline{C}, \overline{C}]$, be a filtering line intersecting $\Gamma(\p) \cup \Gamma(\psi)$ in at least two different points.
    If there exists $P$ in $r_{(a,b)}\cap \alpha_P\cap \beta_P$, where $\alpha_P, \beta_P\in \text{Ctr}(\p, \psi)$ and $\alpha_P\neq\beta_P$ (they might be two overlapping contours), then $(a,b)$ is a special value.
\end{rem}

We recall \cite[Proposition 5.2]{slope1}, showing that the special set is closed. 
This will be fundamental in the proof of our main theorem in the next section.

\begin{prop}\label{spe_closed}
    $\spe{\p}{\psi}$ is a closed subset in $]0,1[ \times [-\overline{C},\overline{C}]$.
\end{prop}

\begin{defn}
The \textbf{ultraspecial set} of $\p$ and $\psi$, denoted by $\mathrm{USp}(\p,\psi)$, is the collection of all $(a,b)$ in $]0,1[ \times [-\overline{C},\overline{C}]$ for which three pair-wise distinct pairs of contours $\{\alpha_P, \alpha_Q\}$,$\{\alpha_R, \alpha_S\}$, $\{\alpha_T, \alpha_U\}$ in $\mathrm{Ctr(\p,\s)}$ intersecting $r_{(a,b)}$ exist, such that every two of them are $(a,b)$-special; that is
    \begin{itemize}
        \item $c_1\lvert x_P-x_Q\rvert=c_2\lvert x_{R}-x_{S}\rvert = c_3\lvert x_{T}-x_{U}\rvert$, with $c_1,c_2,c_3\in \{1,2\}$, if $a\le \frac{1}{2}$,
        \item $c_1\lvert y_P-y_Q\rvert=c_2\lvert y_{R}-y_{S}\rvert=c_3\lvert y_{T}-y_{U}\rvert$, with $c_1,c_2,c_3\in \{1,2\}$, if $a\ge \frac{1}{2}$,
    \end{itemize}
    where $P=P_{(a,b)}=r_{(a,b)}\cap \alpha_P$,
    $Q=Q_{(a,b)}=r_{(a,b)}\cap \alpha_Q$, $R=R_{(a,b)}=r_{(a,b)}\cap \alpha_R$,
    $S=S_{(a,b)}=r_{(a,b)}\cap \alpha_S$, $T=T_{(a,b)}=r_{(a,b)}\cap \alpha_T$, and
    $U=U_{(a,b)}=r_{(a,b)}\cap \alpha_U$ and $x_{*}$, $y_{*}$ denote abscissas
    and ordinates of these points, respectively.
    An element of the ultraspecial set $\mathrm{USp}(\p,\psi)$ is called an \textbf{ultraspecial value} of the pair $(\p, \psi)$.
Any three pairs of contours for which one of the conditions above is satisfied is said $(a,b)$\textbf{-ultraspecial}.
\end{defn}
By definition, $\mathrm{USp}(\p,\psi) \subset \spe{\p}{\psi}$.

\begin{ex}
    Let $\alpha_1, \alpha_2, \alpha_3$ be three disjoint, parallel improper contours in $\textup{Ctr}(\p, \psi)$ such that $\mathrm{dist}(\alpha_1, \alpha_2) = \mathrm{dist}(\alpha_2, \alpha_3)$, where $\mathrm{dist}$ denotes the infimum of the distance between any two points of the contours. Then any $(a,b)$ such that the filtering line $r_{(a,b)}$ meets all three contours is ultraspecial. Indeed, if $P_i = r_{(a,b)} \cap \alpha_i$, $i=1,2,3$, then
    \begin{align*}
        \vert x_{P_1} - x_{P_3} \vert = 2 \vert x_{P_1} &- x_{P_2} \vert = 2 \vert x_{P_2} - x_{P_3} \vert, \\
        \vert y_{P_1} - y_{P_3} \vert = 2 \vert y_{P_1} &- y_{P_2} \vert = 2 \vert y_{P_2} - y_{P_3} \vert,
    \end{align*}
    so $\{\alpha_1,\alpha_2\}, \{\alpha_2,\alpha_3\}$ and $\{\alpha_3,\alpha_1\}$ are $(a,b)$-ultraspecial both when $a \leq \frac{1}{2}$ and $a \geq \frac{1}{2}$.
    In this example the subset of ultraspecial values of $(\p, \s)$  has positive measure.
However, we observe that it is not a generic example.
\end{ex}


\section{Our main result}
\label{sec:mainresult}

In this section we characterise the parameter values at which the matching distance is attained. 
To do this, we define a set $\mathcal{U}(\varphi, \psi)$ containing the ultraspecial set of $\p$ and $\psi$ and another subset of special values.
Explicitly,
\begin{equation}\label{ajfkrnosldwlxl}
\mathcal{U}(\p, \s) = \textup{USp}(\p, \s) \cup (\mathcal{C}\cap \textup{Sp}(\p, \psi))
\end{equation}
where $\mathcal{C}$ is defined in Appendix B. 
Despite the definition of $\mathcal{C}$ plays a fundamental role in the proof of our main result, its definition is rather involved, so we leave its analysis to the appendix.
It is worth mentioning, however, that $\mathcal{C}$ can be approximated by a finite collection of curves, possibly with boundary.

From now on we assume the following to hold:
\begin{ass}\label{uewovneqkdjan}
$\spe{\p}{\psi}$ is a finite union of smooth curves.
\end{ass}

Now we state our main theorem, which we shall prove at the end of the section.

\begin{thm}\label{main}
    Let $\p, \psi \colon M \to \R^2$ be smooth functions satisfying Assumptions \ref{ewaiowjnvornsje} and \ref{uewovneqkdjan}. 
    Then 
    $$
    \Dmatch{\p}{\psi} = \max 
    \left\{\dmatch{\dgm{\p}{a}{b}}{\dgm{\psi}{a}{b}} \mid (a,b) \in \left(\left\{\frac{1}{2}\right\} \times [-\overline C,\overline C]\right) \cup \mathcal{U}(\p,\psi)\right\}.
    $$
\end{thm}

Theorem \ref{main} says that the matching distance is either realised at a parameter value associated with a line of slope 1, or in the set $\mathcal{U}(\p, \psi)$. 
Note that this is an improvement with respect to \cite[Theorem 5.4]{slope1} because $\mathcal{U}(\p, \psi)$ is included in $\text{Sp}(\p, \psi)$.

The proof of Theorem \ref{main} goes as follows. 
We assume by contradiction that the matching distance is realised by filtering lines parameterised by an element outside of $(\{\frac{1}{2}\}\times [-\overline{C}, \overline{C}])\cup \mathcal{U}(\p, \psi)$. 
Among those, we consider the closest $(\overline a, \overline b)$ to the line $a=\frac{1}{2}$. 
Then we show that there exists a perturbation of $(\overline a,\overline b)$, obtained by rotating the corresponding filtering line in the direction of a line of slope 1, that increases the bottleneck distance, contradicting the minimality assumption.
The rest of the section is devoted at showing the details of this procedure. 

Given a line $r_{(a,b)}$, any line $r_{(a',b')}$ can be obtained from $r_{(a,b)}$ either by rotating around a point or by translating it. 
In particular, if $b=b'=0$, the rotation is around the origin, and if $a=a'$, then a translation is being performed.
Thus, the following notation is justified.
\begin{defn}
    For any $(a,b) \ne (a',b') \in [0,1]\times \R$, the symbol $(a,b) \curvearrowright (a',b')$ represents the \textbf{rotation} taking the line $r_{(a,b)}$ to the line $r_{(a',b')}$.
    We say $(a,b) \curvearrowright (a',b')$ is \textbf{clockwise} if $a<a'$, and \textbf{counter-clockwise} if $a'<a$.
    If $a = a'$, $(a,b) \curvearrowright (a',b')$ represents a \textbf{translation}.
\end{defn}

\begin{rem}
Observe that any rotation or translation corresponds to a closed segment $[(a,b),(a',b')] \subset [0,1]\times\R$.
\end{rem}

We recall from Remark~\ref{uwuhhddijdsfjdsfhi} that, if $\dmatch{\dgm{\p}{a}{b}}{\dgm{\s}{a}{b}}>0$, there exist $A$ and $B$ in $r_{(a,b)} \hcap (\Gamma(\p)\cup(\Gamma(\psi))$ and $c \in \{\frac{1}{2},1\}$ such that one of the following holds:
\begin{align*}
& \dmatch{\dgm{\p}{a}{b}}{\dgm{\s}{a}{b}} = c \vert x_A-x_B \vert, \text{ if } a\le \frac{1}{2},\\
& \dmatch{\dgm{\p}{a}{b}}{\dgm{\s}{a}{b}} = c \vert y_A-y_B \vert, \text{ if } a\ge \frac{1}{2}.
\end{align*}


\begin{lem}\label{ver_vert}
For every $b\in  \R$ 
there are $b',b''$ with $b'<b<b''$ such that \begin{align*}
    \dmatch{\dgm{\p}{0}{b}}{\dgm{\psi}{0}{b}} &\leq \dmatch{\dgm{\p}{0}{b'}}{\dgm{\psi}{0}{b'}}, \\
    \dmatch{\dgm{\p}{1}{b}}{\dgm{\psi}{1}{b}} &\leq \dmatch{\dgm{\p}{1}{b''}}{\dgm{\psi}{1}{b''}}.
\end{align*}
\end{lem}

\begin{proof}

We will prove the first inequality since the proof of the second one is analogous.
If $\dmatch{\dgm{\p}{0}{b}}{\dgm{\psi}{0}{b}}=0$ the result is clear. 
Let us assume  that $\dmatch{\dgm{\p}{0}{b}}{\dgm{\psi}{0}{b}}>0$ and, by contradiction, that the statement in the lemma is false, then 
there is a strictly increasing sequence $(b_n)_{n\in \N}$ with $\lim_{n\to\infty}b_n=b$, for which
\begin{equation}\label{eqabs1}
    \dmatch{\dgm{\p}{0}{b_n}}{\dgm{\psi}{0}{b_n}} < \dmatch{\dgm{\p}{0}{b}}{\dgm{\psi}{0}{b}}
\end{equation}
for every $n$.
In~\cite[Theorem 4.3]{slope1}, it was shown that the functions $\p^\ast_{(\cdot, \cdot)}$ and $\psi^\ast_{(\cdot, \cdot)}$ are locally Lipschitz, hence the sequences $(\varphi^*_{(0,b_n)})_{n\in\N}$ and
$(\psi^*_{(0,b_n)})_{n\in\N}$ uniformly converge to the functions $\varphi^*_{(0,b)}$ and
$\psi^*_{(0,b)}$, respectively.
The stability of persistence diagrams allows us to assume that, up to subsequences,
$\dmatch{\dgm{\p}{0}{b_n}}{\dgm{\psi}{0}{b_n}}>0$
for every $n$.

Now, for each $n \in \N$, the generalized Position Theorem (see Remark~\ref{uwuhhddijdsfjdsfhi}) guarantees the existence of two points $A_n,B_n\in r_{(0,b_n)} \hat \cap (\Gamma(\p) \cup \Gamma(\psi))$ such that $$c_n\lvert x_{A_n}-x_{B_n}\lvert=\dmatch{\dgm{\p}{0}{b_n}}{\dgm{\psi}{0}{b_n}}>0$$ with $c_n \in \left\{\frac{1}{2},1\right\}$, for every index $n$.
It is not restrictive to assume that $c_n=c$ and $x_{A_n}< x_{B_n}$ for every $n$.
We observe that, due to the definition of the intersection operator $\hcap$ (Definition~\ref{def:intersection_operator}), while all the points in the real line $r_{(0,b_n)}$ have abscissa $b_n$, it can happen that $x_{A_n},x_{B_n}>b_n$.
Since the set $\Gamma(\p) \cup \Gamma(\psi)$ is compact, it is also not restrictive to assume that two points $\bar A,\bar B$ exist in $\Gamma(\p) \cup \Gamma(\psi)$, such that
$\lim_{n\to\infty} A_n=\bar A$ and $\lim_{n\to\infty} B_n=\bar B$, possibly up to taking subsequences. 

Since all the lines $r_{(0,b_n)}$ are vertical and $b_n\le x_{A_n}< x_{B_n}$, there is a vertical improper contour $v_n\in \mathrm{Ctr}(\p,\psi)$ of abscissa $x_{v_n}$ such that $B_n=(x_{v_n},\infty)\in v_n$ for any $n$. 
The finiteness of the set of improper vertical contours in $\mathrm{Ctr}(\varphi,\psi)$ allows us to assume, up to subsequences, that 
all improper vertical contours $v_n$ coincide with a unique improper contour $v\in \mathrm{Ctr}(\p,\psi)$, having abscissa $x_v$, and that 
$B_n=\bar B=(x_v,\infty)\in v$ for any $n$.

Once again because of the finiteness of the set $\mathrm{Ctr}(\varphi,\psi)$, we can assume, up to subsequences, that one of these cases occurs:
\begin{enumerate}
  \item $A_n\in r_{(0,b_n)}$ and hence $x_{A_n}=b_n$, for any $n$;
  \item there is a vertical improper contour $v'\in \mathrm{Ctr}(\varphi,\psi)$ of abscissa $x_{v'}$ such that $x_{A_n}=x_{v'}$ for every $n$. 
\end{enumerate}


Let us prove that both these cases lead to a contradiction. 

In case 1, $x_{\bar A}=\lim_{n\to\infty} x_{A_n}=\lim_{n\to\infty} b_n=b$ and 
$x_{\bar B}=x_{B_n}=x_v$ for any $n$.
We also recall that $b_n<b$ for any $n$, and observe that $b\le x_v$ (since $b_n=x_{A_n}<x_{B_n}=x_v$ for any $n$, implying that $b=\lim_{n\to\infty} b_n\le x_v$). It follows that 
\begin{align*}
\dmatch{\dgm{\p}{0}{b_n}}{\dgm{\psi}{0}{b_n}}
   &= c\left(x_{B_n}-x_{A_n}\right) \\
   &= c\left(x_{v}-b_n\right) \\
   &> c\left(x_{v}-b\right) \\
   &= \lim_{n\to\infty} c\left(x_v-b_n\right)\\
   &= \lim_{n\to\infty} c\left(x_{B_n}-x_{A_n}\right) \\
   &= \lim_{n\to\infty} \dmatch{\dgm{\p}{0}{b_n}}{\dgm{\psi}{0}{b_n}}\\
   &= \dmatch{\dgm{\p}{0}{b}}{\dgm{\psi}{0}{b}}
\end{align*}
against (\ref{eqabs1}).

In case 2, we have that, for any $n$,
\begin{align*}
\dmatch{\dgm{\p}{0}{b_n}}{\dgm{\psi}{0}{b_n}}
   &= c\lvert x_{A_n}-x_{B_n}\lvert \\
   &= c\lvert x_{v'}-x_v\lvert .
\end{align*}

Because of the stability of persistence diagrams (Theorem \ref{thm:bottleneck_stability})
$$\dmatch{\dgm{\p}{0}{b}}{\dgm{\psi}{0}{b}}=\lim_{n\to\infty} \dmatch{\dgm{\p}{0}{b_n}}{\dgm{\psi}{0}{b_n}}=c\lvert x_{v'}-x_v\lvert $$
and hence, for any $n$,
$$\dmatch{\dgm{\p}{0}{b_n}}{\dgm{\psi}{0}{b_n}}=\dmatch{\dgm{\p}{0}{b}}{\dgm{\psi}{0}{b}}$$
once again against (\ref{eqabs1}).
\end{proof}

\begin{cor} \label{swcoinwcnwkl}
    The suprema of 
    $\dmatch{\dgm{\p}{a}{b}}{\dgm{\psi}{a}{b}}$ for $(a,b)\in [0,1]\times \R$ and for $(a,b)\in ]0,1[\times \R$ concide.
\end{cor}

\begin{proof}
Proposition 4.4 in~\cite{slope1} imply that Lemma \ref{ver_vert} implies that, for any $b> -\overline C$,
$$\dmatch{\dgm{\p}{0}{b}}{\dgm{\psi}{0}{b}}\le \dmatch{\dgm{\p}{0}{-\overline C}}{\dgm{\psi}{0}{-\overline C}}.$$ 
In \cite{slope1} it was proven that, for any $b\le -\overline C$, $$\dmatch{\dgm{\p}{0}{b}}{\dgm{\psi}{0}{b}}= \dmatch{\dgm{\p}{\frac{1}{2}}{b}}{\dgm{\psi}{\frac{1}{2}}{b}}.$$ 

Therefore, for any $b\in \R$,
$$\dmatch{\dgm{\p}{\frac{1}{2}}{-\overline C}}{\dgm{\psi}{\frac{1}{2}}{-\overline C}}\ge \dmatch{\dgm{\p}{0}{b}}{\dgm{\psi}{0}{b}}.$$

We can prove analogously that, for any $b\in \R$,
$$\dmatch{\dgm{\p}{\frac{1}{2}}{\overline C}}{\dgm{\psi}{\frac{1}{2}}{\overline C}}\ge \dmatch{\dgm{\p}{1}{b}}{\dgm{\psi}{1}{b}}.$$

\end{proof}

We now show that we can find a rotation of the filtering line associated with $(a,b)$ that increases the bottleneck distance, if $(a,b)\not\in \spe{\p}{\s}$ and $a\neq 0,\frac{1}{2},1$.
\begin{lem}\label{main_caso1}
Consider $(a,b) \in ]0,1[ \times [-\overline C,\overline C]\setminus \spe{\p}{\psi} $ with $a\neq \frac{1}{2}$.
Then, if $a<\frac{1}{2}$ (resp. $a>\frac{1}{2}$), there exists a clockwise (resp. counter-clockwise) rotation $(a,b) \curvearrowright (a',b')$
such that
\[
    \dmatch{\dgm{\p}{a}{b}}{\dgm{\psi}{a}{b}} \leq \dmatch{\dgm{\p}{a'}{b'}}{\dgm{\psi}{a'}{b'}}.
\]
\end{lem}
\begin{proof}
If $\dmatch{\dgm{\p}{a}{b}}{\dgm{\psi}{a}{b}}=0$ the statement is clear, so we can assume $\dmatch{\dgm{\p}{a}{b}}{\dgm{\psi}{a}{b}}>0$.
Assume, further, that $a< \frac{1}{2}$. Since $(a,b)$ is not a special value, there is a unique pair of points $A,B$ in $r_{(a,b)} \cap (\Gamma(\p) \cup \Gamma(\psi))$ (see Remark~\ref{uwuhhddijdsfjdsfhi}) for which $\dmatch{\dgm{\p}{a}{b}}{\dgm{\psi}{a}{b}}=c\lvert x_A-x_B\rvert$, with $c \in \{\frac{1}{2},1\}$. 
Consider a clockwise rotation $(a,b) \curvearrowright (a',b')$ around the point $A$ such that the segment $[(a,b),(a',b')]$ does not intersect $\spe{\p}{\psi}$.
Such a rotation exists because of Proposition~\ref{spe_closed}.
Thus, the rotation $(a,b) \curvearrowright (a',b')$ induces a bijection
\begin{align*}
    \Psi \colon r_{(a,b)} \cap (\Gamma(\p) \cup \Gamma(\psi)) &\longrightarrow r_{(a',b')} \cap (\Gamma(\p) \cup \Gamma(\psi)) \\
   r_{(a,b)} \cap \alpha_X = X &\longrightarrow \Psi(X) =  r_{(a',b')} \cap \alpha_X.
\end{align*}
This bijection fixes $A = \Psi(A)$ and sends $B$ to $\Psi(B)$. 
Thus, property (iv) in Assumption~\ref{ewaiowjnvornsje} ensures that $\vert x_{\Psi(B)} - x_A \vert \geq \vert x_B - x_A \vert$.
Since there are no special values in the segment $[(a,b),(a',b')]$, one can see that $ \dmatch{\dgm{\p}{a'}{b'}}{\dgm{\psi}{a'}{b'}} = c \vert x_{B'} - x_A \vert$. 
Hence,
\begin{align*}
\dmatch{\dgm{\p}{a'}{b'}}{\dgm{\psi}{a'}{b'}} &=  c \vert x_{B'} - x_A \vert \\
& \geq  c\vert x_B - x_A \vert \\
& = \dmatch{\dgm{\p}{a}{b}}{\dgm{\psi}{a}{b}}.
\end{align*}

The case $a > \frac{1}{2}$ is obtained analogously if replacing abscissas with ordinates and clockwise rotations with counter-clockwise rotations. 
\end{proof}

The next result allows us to find a rotation that increases the bottleneck distance even when the line is associated with special values that are not in $\mathcal{U}(\varphi, \psi)$. 
The proof of this can be found in Appendix B because it depends on the explicit description of the set $\mathcal{C}$.
\begin{lem}\label{main_caso2_1}
Consider $(a,b) \in \spe{\p}{\s} \setminus \mathcal U(\p, \s)$.
There exists a rotation $(a,b)\curvearrowright (a',b')$ for which
\[
\dmatch{\dgm{\p}{a}{b}}{\dgm{\psi}{a}{b}} < \dmatch{\dgm{\p}{a'}{b'}}{\dgm{\psi}{a'}{b'}}.
\]
\end{lem}

We are now ready to prove our main result.

\begin{proof}[Proof of Theorem~\ref{main}]
By contradiction, let $(\overline a, \overline b) \notin \mathcal{U}(\p,\s)$ with $\overline a\neq\frac{1}{2}$ and such that $\Dmatch{\p}{\psi} = \dmatch{\dgm{\p}{\overline a}{\overline b}}{\dgm{\psi}{\overline a}{\overline b}}$.
Among all such $(\overline a, \overline b)$ consider one minimizing $\lvert a-\frac{1}{2}\rvert$. 
Such minimum exists because the collection of parameter values realising the matching distance is compact (see \cite[Section 5]{slope1}).
We assume that $\overline a < \frac{1}{2}$.
By Corollary~\ref{swcoinwcnwkl}, we can also assume that $\overline a>0$. 
Then the following two cases may occur:
\begin{enumerate}

\item $(\overline a,\overline b) \notin \spe{\p}{\psi}$. 
By Lemma~\ref{main_caso1} there is a clockwise rotation $(\overline a, \overline b) \curvearrowright (a',b')$ such that
\[
  \dmatch{\dgm{\p}{\overline a}{\overline b}}{\dgm{\psi}{\overline a}{\overline b}} \le \dmatch{\dgm{\p}{a'}{b'}}{\dgm{\psi}{a'}{b'}}.
\]

Since $(\overline a, \overline b) \curvearrowright (a',b')$ is clockwise, $a'>\overline a$, and so $\vert a' - \frac{1}{2} \vert < \vert \overline a - \frac{1}{2} \vert$, contradicting the minimality of $\overline a$.       
\item $(\overline a,\overline b) \in \spe{\p}{\psi} \setminus \mathcal{U}(\p,\psi)$.
By Lemma \ref{main_caso2_1}, there is a (not necessarily clockwise nor counter-clockwise) rotation $(\overline a, \overline b) \curvearrowright (a',b')$ such that 
\[
        \dmatch{\dgm{\p}{\overline a}{\overline b}}{\dgm{\psi}{\overline a}{\overline b}}<\dmatch{\dgm{\p}{a'}{b'}}{\dgm{\psi}{a'}{b'}} 
\]
contradicting the fact that $(\overline a, \overline b)$ realises the matching distance. 
        
\end{enumerate}
The proof for $\overline a > \frac{1}{2}$ is analogous. 
\end{proof}

\section{Conclusions and future work}

A central problem in the study of $2$-parameter persistence is the reduction of the computational cost of the matching distance. 
In this article we restrict the set of filtering lines required for its computation.
We pass from a rectangular parameter set to the union of one closed segment, corresponding to lines of slope 1, and a subset $\mathcal{U}(\p,\s)$ providing an obstruction for the uniqueness of an optimal matching when computing the aforementioned metric.
We conjecture the set $\mathcal{U}(\p,\s)$ to be finite in the class of regular filtering functions.
This conjecture can be split into two parts. 
First, we conjecture that $\textup{USp}(\p, \psi)$ is finite. 
Second, if both $\textup{Sp}(\p, \psi)$ and $\mathcal{C}$ are finite collections of curves, then their intersection is generically finite. 
Moreover, it remains open the problem of exploiting the characterisation of the matching distance in Theorem~\ref{main} to design a new algorithm to compute it.

\bibliographystyle{plainurl}
\bibliography{biblio_HMD}

\section{Appendix A}\label{appendixA}
The following results about persistence diagrams are well-known (see, for example, \cite{pers_Betti_stable, natural_opt_matching}). 
However, we provide the proofs for self-containment.

\begin{lem}\label{loc_finiteness}
Let $\p\colon X\to \R$ be a continuous function from a compact topological space $X$. Then
\begin{itemize}
\item $\mathrm{Dgm}(\p)$ contains a finite number of essential cornerpoints, 
\item for every finite $\eps>0$, $\mathrm{Dgm}(\p)\setminus U_\eps$ contains a finite number of proper cornerpoints, where $U_\eps$ is an $\varepsilon$-neighborhood of $\Delta$ in $(\semip,d)$.
\end{itemize}
\end{lem}
\begin{proof}
By contradiction, let us assume that
$\Dgm{\p}$ contains an infinite set $S_\infty$ of essential cornerpoints.
The compactness of $X$ and the definition of multiplicity $m(u,v)$ ensure us that if $(u,\infty)\in S_\infty$, then  $\min \p\le u\le \max \p$.
Therefore, the set $C$ of the abscissas of points in $S_\infty$ admits an accumulation point in the interval $[\min\p,\max\p]$.
Since the PBNF $\beta^\p_k(u,v)$ is non-decreasing in $u$ and non-increasing in $v$,
the function $\beta_p^{\p}$ takes the value $\infty$ at some point in $\Delta^+$, against Assumption~\ref{jdshdsds}.

Let us show the second statement by contradiction. Assume that a real $\varepsilon >0$ exists such that there is an infinite subset $S$ in $\mathrm{Dgm}(\p)\setminus U_\eps$ of proper cornerpoints
By the compactness of $X$ and the definition of multiplicity $m(u,v)$, every $(u,v) \in S$ must be such that $\min \p \le u < v \le \max \p$. 
Therefore $S \subset K = \{(u,v) \in \semip \colon \min \p \le u \le v \le \max \p\}$.
By compactness of $K$, there exists an accumulation point $(u',v')$ of $S$ in $K$.
Therefore, for any $0<\delta<\varepsilon$ the open box $]u'-\delta, u'+\delta[ \times ]v'-\delta, v'+\delta[$ centered in $(u',v')$ contains an infinite number of proper cornerpoints.
Since the PBNF $\beta^\p_k(u,v)$ is non-decreasing in $u$ and non-increasing in $v$, there exists $(u'',v'') \in ]u'-\delta, u'+\delta[ \times ]v'-\delta, v'+\delta[$ such that $\beta^\p_k(u'',v'')=\infty$, contradicting Assumption \ref{jdshdsds}.
\end{proof}

As an immediate consequence, we have that the points of positive multiplicity of a persistence diagram, often referred as the persistence diagram itself, form a compact set.
\begin{prop}\label{cpt}
Let $\p\colon X\to \R$ be a continuous function from a compact topological space $X$. 
Then $\mathrm{Dgm}(\p)$ is a compact subset of $(\semip, d)$.
\end{prop}

Under the assumptions of this article, an optimal matching for the bottleneck distance always exists (confront with~\cite[Theorem 28]{natural_opt_matching}).
\begin{prop}\label{opt_matching}
For every $\p, \psi \colon X \to \R$ continuous functions from a compact topological space $X$
\[
\dmatch{\textup{Dgm}(\p)}{\textup{Dgm}(\psi)} = \min_{\sigma}\max_{p\in \mathrm{Dgm}(\p)}d(p,\sigma(p)).
\]
\end{prop}
\begin{proof}
By Proposition~\ref{cpt}, for every bijection $\sigma$, $\text{cost}(\sigma)=\max_{p\in \Dgm{\p}}d(p,\sigma(p))$.
Thus, assume $\dmatch{\Dgm{\p}}{\Dgm{\s}} = D<\infty$. 
By contradiction suppose cost$(\sigma) > D$ for every matching $\sigma \colon \Dgm{\p} \to \Dgm{\s}$.
Then there is a sequence  of matchings $(\sigma_n)_{n \in \N}$ such that $0 < \textup{cost}(\sigma_n) - D < \frac{1}{n}$.
Let $U_D$ be the open ball of radius $D$ centered at $\Delta$ in $(\bar \Delta^*, d)$. 
For each bijection $\sigma$ consider the injection $\sigma_1\colon \Dgm{\p}\to \Dgm{\s}$ sending 
$p$ to $ \Delta$, if $p \in U_D \cap \Dgm{\p}$, and $
p$ to $ \sigma(p)$, otherwise.
In this way cost$(\sigma_1) \le \textup{cost}(\sigma)$, for each $\sigma$.
Analogously, consider the injection $\sigma_2\colon  \Dgm{\s}\to  \Dgm{\p}$ sending $\sigma(p)$ to $\Delta$, if $\sigma(p)\in U_D\cap \Dgm{\s}$, and $\sigma(p)$ to $p$, otherwise. 
We have also in this case that cost$(\sigma_2) \le \textup{cost}(\sigma)$, for each $\sigma$.
By the proof of Cantor-Bernstein-Schröder theorem, we can construct a bijection $\overline{\sigma}\colon  \Dgm{\p}\to  \Dgm{\s}$ such that 
cost$(\overline{\sigma})$ is
smaller or equal to cost$(\sigma)$. 
Observe, in addition, that if cost$(\overline \sigma) = d(p,\Delta)$, for some $p \in U_D \cap  \Dgm{\p}$, then cost$(\overline \sigma) < D$, against the assumption that $D = \inf_{\sigma} \textup{cost}(\sigma)$.
Thus, for each $\overline \sigma$, cost$(\overline \sigma)$ belong to the set $\{d(p, \sigma(p))\}_{p \in  \Dgm{\p} \setminus U_D}$.
But this set is finite because of Lemma \ref{loc_finiteness}.
In particular, this implies that the cardinality of the set $ \{\textup{cost}( \sigma_n)\}_{n \in \N}$ is finite.
Then there exists $N \in \N$ sufficiently big such that, for each $n > N$, $\textup{cost}(\sigma_n) = D$, which is a contradiction. 

\end{proof}

\section{Appendix B}\label{appendixB}
Our goal is now to show the following:
\begin{thm*}
Consider $(a,b) \in \spe{\p}{\s} \setminus \mathcal U(\p, \s)$.
There exists a rotation $(a,b)\curvearrowright (a',b')$ for which
\[
\dmatch{\dgm{\p}{a}{b}}{\dgm{\psi}{a}{b}} < \dmatch{\dgm{\p}{a'}{b'}}{\dgm{\psi}{a'}{b'}}.
\]
\end{thm*}

To do this, we will exhibit an explicit approximation of the set $\mathcal{C}\subset \mathcal{U}(\p, \psi)$ in Section~\ref{sec:mainresult}. 
What follows is devoted to this goal.

Consider 
four contours $\alpha_P, \alpha_Q, \alpha_R, \alpha_S \in \textup{Ctr}(\p, \s)$ such that $\{\alpha_P,\alpha_Q\} \ne \{\alpha_R,\alpha_S\}$.
Define the set $\mathcal{Q}$ of $(a,b) \in \left]0,\frac{1}{2}\right]\times[-\overline{C},\overline{C}]$ such that the line $r_{(a,b)}$ intersects the contours $\alpha_P,\alpha_Q,\alpha_R,\alpha_S$, respectively, in $P(a,b)$, $Q(a,b)$, $R(a,b)$, $S(a,b)$. 
Note that such intersections, if they exist, are unique.
Define the functions
\begin{align}\label{fandg}
\nonumber
     f=f_{\alpha_P,\alpha_Q} \colon \mathcal{Q}&\to [0,\infty[\\ 
    (a,b) &\mapsto (x_{P(a,b)}-x_{Q(a,b)})^2, \\ \nonumber
     g=g_{\alpha_R,\alpha_S} \colon \mathcal{Q} &\to [0,\infty[\\ \nonumber
    (a,b) &\mapsto (x_{R(a,b)}-x_{S(a,b)})^2.
\end{align}

Similarly, consider 
four contours $\beta_{P'}, \beta_{Q'}, \beta_{R'}, \beta_{S'} \in \textup{Ctr}(\p, \s)$ such that $\{\beta_{P'},\beta_{Q'}\} \ne \{\beta_{R'},\beta_{S'}\}$.
Define the set $\mathcal{Q'}$ of $(a,b) \in \left[\frac{1}{2}, 1\right[\times[-\overline{C},\overline{C}]$ such that the line $r_{(a,b)}$ intersects the contours $\beta_{P'}, \beta_{Q'}, \beta_{R'}, \beta_{S'}$, respectively, in $P'(a,b)$, $Q'(a,b)$, $R'(a,b)$, $S'(a,b)$. 
Note that such intersections, if they exist, are unique.
Define the functions
\begin{align}\label{fandg2}
\nonumber
     f'=f'_{\beta_{P'},\beta_{Q'}} \colon \mathcal{Q}'&\to [0,\infty[\\ 
    (a,b) &\mapsto (y_{P'(a,b)}-y_{Q'(a,b)})^2, \\ \nonumber
     g'=g'_{\beta_{R'},\beta_{S'}} \colon \mathcal{Q}' &\to [0,\infty[\\ \nonumber
    (a,b) &\mapsto (y_{R'(a,b)}-y_{S'(a,b)})^2.
\end{align}

\begin{rem}
The functions $f$ and $g$ measure the square of the distances of the abscissas and of the ordinates of the pairs of points $P(a,b), Q(a,b)$ and $R(a,b), S(a,b)$ when $(a,b)$ varies. 
It is possible to show that $f$ and $g$ are $C^1$ on the interior of $\mathcal Q$.
Analogously, $f'$ and $g'$ are also $C^1$ on the interior of their domain.
\end{rem}

Let us define the set
\begin{align*}
\mathcal{C} & = \left\{ (a, b) \in \left]0,\frac{1}{2}\right]\times \R \mid
\begin{array}{l}
\nabla f, \nabla g \textup{ are parallel or } r_{(a,b)}  \\
\textup{ intersects a contour at an endpoint}
\end{array}
\right\}\\
&\cup \left\{ (a, b) \in \left[\frac{1}{2}, 1\right[\times \R \mid
\begin{array}{l}
\nabla f', \nabla g' \textup{ are parallel or } r_{(a,b)}\\
\textup{ intersects a contour at an endpoint}
\end{array}
\right\}
\end{align*}
\textit{Claim:} $\mathcal{C}$ can be approximated by a finite union of curves, possibly with boundary.

First, observe that the collection of $(a,b)$ associated with lines through the endpoint of a contour is a segment in the parameter space. 
By Assumption~\ref{ewaiowjnvornsje}, the number of contours is finite, so there is just a finite number of such lines.
Let us now study when $\nabla f =\lambda \nabla g$ and $\nabla g= \mu \nabla f$ for some $\lambda$ and $\mu$ real numbers. 
This happens when 
$\frac{\partial f}{\partial a}\frac{\partial g}{\partial b}-\frac{\partial f}{\partial b}\frac{\partial g}{\partial a}=0$. 
We consider an approximation of this equation obtained by calculating the abscissas of the intersection between a filtering line $r_{(a,b)}$ given by $y = \frac{1-a}{a}x - \frac{b}{a}$ and the line $y=m_\omega (x-x_X)+y_X$ tangent to the contour $\omega$ in $\{\alpha_P, \alpha_Q, \alpha_R, \alpha_S\}$ at $X \in r_{(a,b)} \cap \omega$, for each $(a,b)$ in $\mathcal Q$.
This provides an approximation of the first order of $r_{(a,b)} \cap \omega$.


Substituting the approximation of the first order of the abscissas of the intersection points on the equation $\frac{\partial f}{\partial a}\frac{\partial g}{\partial b}-\frac{\partial f}{\partial b}\frac{\partial g}{\partial a}=0$ yields the following equation when $m_\omega > - \infty$ for each $\omega \in \{\alpha, \beta, \gamma, \delta\}$, where $\alpha, \beta, \gamma, \delta$ replace, respectively, $\alpha_P, \alpha_Q, \alpha_R, \alpha_S$ for convenience. This can be implemented on any symbolic computation software such as \cite{sagemath}:

\begin{align}\label{cnvjrovnapo}
& \left(\frac{{\left(m_{\alpha} x_{P} - y_{P}\right)} a - b}{a {\left(m_{\alpha} + 1\right)} - 1} - \frac{{\left(m_{\beta} x_{Q} - y_{Q}\right)} a - b}{a {\left(m_{\beta} + 1\right)} - 1}\right)\nonumber\\
&\cdot\left(\frac{{\left(m_{\delta} x_{S} - y_{S}\right)} a - b}{a {\left(m_{\delta} + 1\right)} - 1} - \frac{{\left(m_{\gamma} x_{R} - y_{R}\right)} a - b}{a {\left(m_{\gamma} + 1\right)} - 1}\right) \nonumber\\
&\cdot \left(\frac{m_{\gamma} x_{R} - y_{R}}{a {\left(m_{\gamma} + 1\right)} - 1} - \frac{m_{\delta} x_{S} - y_{S}}{a {\left(m_{\delta} + 1\right)} - 1}\right. \nonumber \\ 
&\left.+\frac{{\left({\left(m_{\delta} x_{S} - y_{S}\right)} a - b\right)} {\left(m_{\delta} + 1\right)}}{{\left(a {\left(m_{\delta} + 1\right)} - 1\right)}^{2}} -\frac{{\left({\left(m_{\gamma} x_{R} - y_{R}\right)} a - b\right)} {\left(m_{\gamma} + 1\right)}}{{\left(a {\left(m_{\gamma} + 1\right)} - 1\right)}^{2}}\right) \nonumber\\
&\cdot\left(\frac{1}{a {\left(m_{\alpha} + 1\right)} - 1} - \frac{1}{a {\left(m_{\beta} + 1\right)} - 1}\right)\nonumber\\
= - &\left(\frac{{\left(m_{\alpha} x_{P} - y_{P}\right)} a - b}{a {\left(m_{\alpha} + 1\right)} - 1} - \frac{{\left(m_{\beta} x_{Q} - y_{Q}\right)} a - b}{a {\left(m_{\beta} + 1\right)} - 1}\right)\\
&\cdot\left(\frac{{\left(m_{\delta} x_{S} - y_{S}\right)} a - b}{a {\left(m_{\delta} + 1\right)} - 1} - \frac{{\left(m_{\gamma} x_{R} - y_{R}\right)} a - b}{a {\left(m_{\gamma} + 1\right)} - 1}\right)\nonumber \\
& \cdot\left(\frac{m_{\alpha} x_{P} - y_{P}}{a {\left(m_{\alpha} + 1\right)} - 1} - \frac{m_{\beta} x_{Q} - y_{Q}}{a {\left(m_{\beta} + 1\right)} - 1} \right. \nonumber \\
& \left.- \frac{{\left({\left(m_{\alpha} x_{P} - y_{P}\right)} a - b\right)} {\left(m_{\alpha} + 1\right)}}{\left(a {\left(m_{\alpha} + 1\right)} - 1\right)}^{2}+ \frac{{\left({\left(m_{\beta} x_{Q} - y_{Q}\right)} a - b\right)} {\left(m_{\beta} + 1\right)}}{{\left(a {\left(m_{\beta} + 1\right)} - 1\right)}^{2}}\right)\nonumber \\
&\cdot\left(\frac{1}{a {\left(m_{\delta} + 1\right)} - 1} - \frac{1}{a {\left(m_{\gamma} + 1\right)} - 1}\right). \nonumber
\end{align}
We make the following observations: 
\begin{enumerate}
    \item When any contour $\omega \in \{\alpha, \beta, \gamma, \delta\}$ is improper vertical, i.e. when $m_\omega = - \infty$, the equation is well defined and it is given by the limit of Equation \ref{cnvjrovnapo} when $m_\omega$ goes to $- \infty$.
    \item For each contour $\omega \in \{\alpha, \beta, \gamma, \delta\}$ and for each $a \in ]0,\frac{1}{2}]$, since the local slope $m_\omega$ is nonpositive, we have $a(m_\omega+1)-1 < 0$.
    Hence the denominators in Equation \ref{cnvjrovnapo} do not vanish and it can be written in the form 
    $$
    C(a,b) D(a,b) (P_1(a,b) P_2(a,b) - Q_1(a,b) Q_2(a,b)) = 0,
    $$
    for $C,D,P_1,P_2,Q_1,Q_2 \in \R[a,b]$.

\item The polynomials $C(a, b)$ and $D(a, b)$ are of the form $c_2a^2 + c_1ab +c_0a$, for some $c_i \in \R$. 
Thus, for each of these terms, the zeros are contained in two lines, one of them being $a=0$, which is not an admissible solution. 
The polynomial $P_1(a,b) P_2(a,b) - Q_1(a,b) Q_2(a,b)$ is zero on a, possibly reducible, algebraic curve. 
The union of these lines and curves is a solution of Equation \ref{cnvjrovnapo}.
\item The procedure above can be replicated by intersecting the filtering line $r_{(a,b)}$ with the Taylor approximation of the contours $\alpha, \beta, \gamma, \delta$ truncated at any finite order to get an arbitrarily precise approximation of the curves in  $\mathcal{C}$.
\end{enumerate}

Now we prove Lemma~\ref{main_caso2_1}:

\begin{proof}[Proof of Lemma 5.8]
Assume that $\dmatch{\dgm{\p}{a}{b}}{\dgm{\psi}{a}{b}}=\lvert x_P-x_Q\rvert=\lvert x_R-x_S\rvert$, where $P,Q,R$ and $S$ belong, respectively, to the $(a,b)$-special contours $\{\alpha_P, \alpha_Q\}, \{\alpha_R,\alpha_S\}$.
Note that since $(a,b)$ is not ultraspecial, 
there is no triple of $(a,b)$-ultraspecial contours.

By hypothesis, $\nabla f$ and $\nabla g$ are nonzero and not parallel. 
Thus, there exists a vector $v=\lambda \nabla f+\mu \nabla g$ with $\lambda, \mu>0$, such that the scalar products $\nabla f\cdot v$ and $ \nabla g\cdot v$ are strictly positive.
By definition of directional derivative, we have that
\[
\frac{\partial f}{\partial v}=\lim_{\varepsilon\to 0}\frac{f((a,b)+\varepsilon v)-f(a,b))}{\varepsilon}>0,
\;\;\;\;
\frac{\partial g}{\partial v}=\lim_{\varepsilon\to 0}\frac{g((a,b)+\varepsilon v)-g(a,b))}{\varepsilon}>0.
\]
In particular, this implies that there exists $\varepsilon _v>0$ such that $f((a,b)+\varepsilon v)>f(a,b)$ and $g((a,b)+\varepsilon v)>g(a,b)$, for every $\varepsilon <\varepsilon_v$. 
Choose one such $\varepsilon$ and define $(a',b')=(a,b)+\varepsilon v $.
By Assumption~\ref{uewovneqkdjan}, $\spe{\p}{\psi}$ is a finite union of curves, hence, $v$ and $\eps < \eps_v$ can be chosen such that $(a',b') \notin \spe{\p}{\psi}$ and $a \ne a'$.
By definition of $f$ and $g$, the filtering line $r_{(a',b')}$ still intersects the contours $\alpha_P, \alpha_Q,\alpha_R,\alpha_S$. 
Since the bottleneck distance is stable (Theorem~\ref{stability}), there are two points $P',Q'$ in the intersections of either $\alpha_P,\alpha_Q$ or $\alpha_R,\alpha_S$ with the line $r_{(a',b')}$ such that $\dmatch{\dgm{\p}{a'}{b'}}{\dgm{\psi}{a'}{b'}}=\lvert x_{P'}-x_{Q'}\rvert$. 
The condition on the positivity of the directional derivatives implies, hence, that $\lvert x_{P'}-x_{Q'}\rvert>\lvert x_{P}-x_{Q}\rvert$. 
Thus, we have found a rotation $(a,b) \curvearrowright (a',b')$ that strictly increases the bottleneck distance:
\[
\dmatch{\dgm{\p}{a}{b}}{\dgm{\psi}{a}{b}} < \dmatch{\dgm{\p}{a'}{b'}}{\dgm{\psi}{a'}{b'}}.
\]

The proof for the case $a \ge \frac{1}{2}$ is the same substituting abscissas with ordinates and the functions $f,g$ with $f',g'$.
\end{proof}

\end{document}